\documentclass{amsart}

\usepackage{amsmath,latexsym,amssymb,amsthm,graphicx}
\usepackage{hyperref}
\usepackage[all]{xy}

\newtheorem{theorem}{Theorem}[section]

\newtheorem{corollary}[theorem]{Corollary}

\newtheorem*{btheorem}{Theorem}
\newtheorem{lemma}[theorem]{Lemma}

\theoremstyle{definition}
\newtheorem{remark}[theorem]{Remark}

\newtheorem{example}[theorem]{Example}

\newcommand{\tn}[1]{\textnormal{#1}}
\newcommand{\wt}[1]{\widetilde{#1}}
\newcommand{\pa}[1]{\left(#1\right)}
\newcommand{\ts}[1]{\left|#1\right|}

\def\CP{{S^{2n-1}/S^{1}}}
\def\Z{\mathbb{Z}}
\def\R{\mathbb{R}}
\def\N{\mathbb{N}}
\def\ss{\vskip 1.5 pt plus 1pt}
\def\M{{M/\hspace{-3pt}\sim}}

\font\cuf=cmtt8
\newcommand{\curl}[1]{{\cuf #1}}

\begin{document}
\title{On fundamental groups of quotient spaces}

\author[J.S.~Calcut]{Jack S. Calcut}
\address{Department of Mathematics\\
         Oberlin College\\
         Oberlin, OH 44074}
\email{jcalcut@oberlin.edu}
\urladdr{\href{http://www.oberlin.edu/faculty/jcalcut/}{\curl{http://www.oberlin.edu/faculty/jcalcut/}}}

\author[R.E.~Gompf]{Robert E. Gompf}
\address{Department of Mathematics\\
                    University of Texas at Austin\hfill\break
                    \indent 1 University Station C1200\\
                    Austin, TX 78712-0257}
\email{gompf@math.utexas.edu}
\urladdr{\href{http://www.ma.utexas.edu/users/gompf/}{\curl{http://www.ma.utexas.edu/users/gompf/}}}

\author[J.D.~McCarthy]{John D. McCarthy}
\address{Department of Mathematics\\
                    Michigan State University\\
                    East Lansing, MI 48824-1027}
\email{mccarthy@math.msu.edu}
\urladdr{\href{http://www.math.msu.edu/~mccarthy/}{\curl{http://www.math.msu.edu/\textasciitilde mccarthy/}}}

\thanks{The second author was partially supported by NSF grants
DMS-0603958 and 1005304.}
\date{September 21, 2011}

\begin{abstract}
In classical covering space theory, a covering map induces an injection of fundamental groups. This paper reveals a dual property for certain quotient maps having connected fibers, with applications to orbit spaces of vector fields and leaf spaces in general.
\end{abstract}

\maketitle

\section{Introduction}\label{s:intro}

Given a map $f:X\to Y$ of topological spaces, an important problem is to relate the homology of $X$ with that of $Y$ using properties of $f$. Vietoris pioneered the study of this problem with his mapping theorem.

\begin{btheorem}[Vietoris~{\cite[\S~III]{vietoris}}]
Let $f:X\to Y$ be a surjective map of compact metric spaces. If the (Vietoris) homology of each fiber $\wt{H}_{r}\pa{f^{-1}\pa{y}}$ is trivial for $0\le r \le n-1$, then the induced homomorphism $f_{\ast}:\wt{H}_{r}\pa{X}\rightarrow\wt{H}_{r}\pa{Y}$ is an isomorphism for $r\le n-1$ and is surjective for $r=n$.
\end{btheorem}

There is an extensive literature on this theorem and its various generalizations\footnote{See, e.g., \cite{leray2}, \cite[pp.~123--124]{dieudonne}, \cite{begle55}, \cite{reitberger02}, \cite[pp.~75, 221, 317]{bredon}, \cite{lean}, \cite{dydak_walsh}, \cite{dydak_kozlowski}, \cite{snyder}, and references therein.}.
Smale first discovered an analog of Vietoris' mapping theorem valid for homotopy groups~\cite{smale_b}. One of Smale's results for the fundamental group is the following.\ss

\begin{btheorem}[Smale~{\cite[Thm.~9,~$n=1$]{smale_b}}]
Let $f:\pa{X,x_{0}}\rightarrow\pa{Y,y_{0}}$ be a proper surjective map of locally compact, locally path-connected, separable metric spaces where $Y$ is semilocally simply-connected. If each fiber $f^{-1}(y)$ is locally path-connected and path-connected, then the induced homomorphism $f_{\sharp}:\pi_{1}\pa{X,x_{0}}\rightarrow\pi_{1}\pa{Y,y_{0}}$ is surjective.
\end{btheorem}

One purpose of this note is to prove the following generalization of Smale's theorem.\ss

\begin{theorem}\label{mainthm}
Let $f:\pa{X,x_{0}}\rightarrow\pa{Y,y_{0}}$ be a quotient map of topological spaces, where $X$ is locally path-connected and $Y$ is semilocally simply-connected. If each fiber $f^{-1}\pa{y}$ is connected, then the induced homomorphism $f_{\sharp}:\pi_{1}\pa{X,x_{0}}\rightarrow\pi_{1}\pa{Y,y_{0}}$ is surjective.
\end{theorem}

\begin{remark}
In Theorem~\ref{mainthm}, no separation axioms are assumed (e.g., Hausdorff) and neither $X$ nor $Y$ is required to be path-connected.
By definition, a space $Z$ is \emph{semilocally simply-connected} if each $z\in Z$ has a neighborhood $U$ such that the induced homomorphism $\pi_{1}\pa{U,z}\to\pi_{1}\pa{Z,z}$ is trivial. We mention that the hypotheses of Theorem~\ref{mainthm} imply that $Y$ is locally path-connected (see Section~\ref{sec:proof}).  Hence, for such $Y$, an alternative notion of semilocal simple connectivity, due to Spanier, coincides with the common definition above~\cite[Thm.~2.8]{frvz}.
\end{remark}

To see that Theorem~\ref{mainthm} implies Smale's quoted theorem, note that a proper map $f:X\to Y$, where $Y$ is locally compact and Hausdorff, is necessarily closed (one may prove this fact using nets), and a closed surjective map is a quotient map. Theorem~\ref{mainthm} also implies parts of \cite[Thm.~2]{kozlowski}, \cite[Thm.~6.1]{dugundji}, and \cite[Thm.~1]{mcmillanshrikhande}.\ss

The second purpose of this note is to present a collection of pertinent examples (see \S2) which we now outline. Theorem~\ref{mainthm} merely requires that fibers be connected, rather than path-connected. Respecting Gromov's dictum that empty generalization should be avoided \cite[p.~339]{berger}, Example~\ref{instances} gives instances of Theorem~\ref{mainthm} with non-path-connected fibers. Examples~\ref{connfibers}--\ref{Xlpc} and~\ref{Yslsc} show that Theorem~\ref{mainthm} becomes false if any hypothesis is omitted. Example~\ref{Yslsc} is probably the simplest imaginable example satisfying all the hypotheses of Theorem~\ref{mainthm} except for semilocal simple-connectivity of the quotient, but where $\pi_1$ is created in the quotient. Remark~\ref{Yslsc_extend} extends Example~\ref{Yslsc} to a similar example where $\pi_1$ is created, but where the domain is simply-connected. Previously, Bing discovered an example which demonstrates this same phenomenon~\cite{bing57}. Bing's example, announced in 1955 (see MR857209), consists of a solenoid $\Sigma\subset\R^3$ and the quotient map $f:\R^3\to\R^3/\Sigma$ crushing $\Sigma$ to a point. Bing's proof that the quotient $\R^3/\Sigma$ is not simply connected nor locally simply-connected appeared later in~\cite{bing86}. Theorem~\ref{mainthm} above now implies that $\R^3/\Sigma$ is not semilocally simply-connected. This may also be deduced from Bing's \emph{proof} of~\cite[Thm.~1]{bing86}. Filling in the details of Bing's argument, however, is nontrivial. In this regard, our conceptually simple (and 2-dimensional) example may be of interest.
For another proof that $\R^3/\Sigma$ is not simply-connected, see~\cite{shrikhande} (cf.~\cite{mcmillan}, \cite{mcmillanshrikhande}). In~\cite{karimov_repovs}, $\R^3/\Sigma$ was shown to have an uncountable first integral homology group (cf. Remark~\ref{Yslsc_extend} below).\ss

For the context of Example~\ref{higherhomotopy}, consider a familiar special case of Theorem~\ref{mainthm}, the projection $p:E\to B$ of a fiber bundle. The associated long exact sequence shows that if the fibers of $p$ are path-connected, then $p_{\sharp}:\pi_{1}(E)\to\pi_{1}(B)$ is surjective. The long exact sequence also gives information about higher homotopy, for example yielding an isomorphism when suitable homotopy of the fiber vanishes. Similarly, Smale's result quoted above extends to higher homotopy groups (cf.~\cite[Thm.~2]{kozlowski} and~\cite[Thm.~6.1]{dugundji}). Such extension implies that the hypotheses of these theorems are necessarily more restrictive than ours. In Example~\ref{higherhomotopy}
we exhibit a quotient map with contractible, locally contractible fibers and
domain, whose quotient is the $2$--sphere and so has higher homotopy (although it is simply connected in accordance with Theorem~\ref{mainthm}).
This shows that in its full generality, Theorem~\ref{mainthm} can only apply to the first homotopy group. In fact, we obtain the 2--sphere as the quotient of a partition of closed upper $3$--space into connected arcs (one open, the rest closed). This example is very simple analytically.\ss

We discovered Theorem~\ref{mainthm} prior to our awareness of the Vietoris mapping theorem, its various extensions, and Bing's example. Our motivation was a question of Arnol{\textprime}d on exotic $\R^4$s appearing as orbit spaces of certain vector fields on $\R^5$ (see Section~\ref{conclusion}). Example~\ref{hopf} below shows that even quadratic polynomial vector fields on euclidean space can have topologically interesting \emph{closed} manifold quotients. Once we observed this example, we sought obstructions to manifolds appearing as orbit spaces of vector fields on a given manifold. In this way, we were led to Theorem~\ref{mainthm}. A subsequent search of the literature revealed a connection to the rich history of Vietoris' mapping theorem, Smale's theorems, and Bing's example. However, the extra generality of Theorem~\ref{mainthm} remains indispensable in this context. For example, quotient maps from $\R^n$ to a compact space are never proper as required by Smale's theorem.\ss

It is not clear whether Theorem~\ref{mainthm} always applies to orbit space quotient maps for arbitrary smooth vector fields. In ~\cite{cg}, the first two authors answer this in the affirmative for generalized gradient vector fields for Morse functions on smooth manifolds $M$. When $\tn{dim}M>0$, these quotient maps are never proper (consider any single nontrivial orbit) and their targets are typically non-Hausdorff, so they also provide additional examples where the generality of Theorem~\ref{mainthm} is useful.\ss

At the level of proofs, Theorem~\ref{mainthm} is not extraordinary. Nevertheless, it appears to have gone unnoticed despite its fundamental nature and the appearance of similar results, e.g., \cite[Cor.~2]{smale_a}.
We mention an obvious duality (for decent spaces): a classical covering map has discrete fibers and induces an injection on $\pi_1$, whereas, by Theorem~\ref{mainthm}, a quotient map with connected fibers induces a surjection on $\pi_1$. Not surprisingly, our proof of Theorem~\ref{mainthm} uses an appropriate cover of $Y$.\ss

This note is organized as follows. Section~\ref{sec:examples} presents examples, Section~\ref{sec:proof} proves Theorem~\ref{mainthm}, Section~\ref{simplicial_category} makes some observations on Theorem~\ref{mainthm} and the simplicial category, and Section~\ref{conclusion} closes with remarks on Arnol{\textprime}d's problem and questions for further study.\ss

By our convention, a \emph{map} is a continuous function, and $A\approx B$ indicates that $A$ is homeomorphic to $B$.\ss

\section{Examples}\label{sec:examples}

This section presents several examples relevant to Theorem~\ref{mainthm}.

\begin{example}\label{connfibers}
The quotient map \[f:[0,1]\to [0,1]/\{0,1\} \approx S^1\] shows that Theorem~\ref{mainthm} minus the hypothesis that fibers are connected is false. For more examples, consider any nontrivial classical covering map.
\end{example}

\begin{example}\label{quomap}
The surjective map $f:[0,1)\to S^1$ given by $f(x)=\exp(2\pi ix)$ shows that Theorem~\ref{mainthm} minus the hypothesis that $f$ is a quotient map is false.
\end{example}

\begin{example}\label{Xlpc}
This example shows that Theorem~\ref{mainthm} minus the hypothesis that $X$ is locally path-connected is false.
Consider the Warsaw circle $W\subset\R^2$ shown in Figure~\ref{warsaw}.
\begin{figure}[h!]
    \centerline{\includegraphics{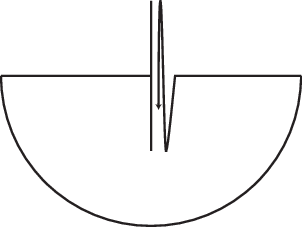}}
    \caption{Warsaw circle $W\subset\mathbb{R}^{2}$.}
    \label{warsaw}
\end{figure}
It equals the disjoint union $A\sqcup B\sqcup C$ where
\begin{align*}
    A&=\left\{(0,y)\mid -1\leq y \leq 1\right\},\\
    B&=\left\{(x,\sin(1/x))\mid 0<x\leq 1/\pi\right\},  \text{ and}\\
    C&=\left\{(x,0)\mid -2\leq x<0\text{ or }1/\pi<x\leq 2\right\}  \cup    \left\{(x,y)\mid    x^{2}+y^{2}=4\text{ \& }y\leq0\right\}.
\end{align*}
Note that $A\cup B$ is connected, $W$ is not locally path-connected, and $W$ is simply-connected. Consider the quotient map
\[
	f:W\to W/(A\cup B) \approx S^1.
\]
\end{example}

\begin{remark}
Example~\ref{Xlpc} also shows that in Theorem~\ref{mainthm} the hypothesis `$X$ is locally path-connected' cannot be replaced with the hypothesis `$Y$ is locally path-connected'.
\end{remark}

\begin{remark}
Example~\ref{Xlpc} may be modified to produce a similar example having all fibers path-connected and not merely connected. Let $\sim$ be the equivalence relation on $W$ with two non-singleton classes, namely $A$ and $B$. Consider the quotient map $f:W\to W/\hspace{-4pt}\sim$. The quotient space is non-Hausdorff, yet it is homotopy equivalent to $S^1$ (see~\cite{cgm} for details).
\end{remark}

\begin{example}\label{Yslsc}
This key example shows that Theorem~\ref{mainthm} minus the hypothesis that $Y$ is semilocally simply-connected is false.
Let $\N$ denote the natural numbers, $\N_{0}=\N\cup\left\{0\right\}$, $I=[0,1]\subset\mathbb{R}$, and $I_{n}=I\times\left\{n\right\}$. Consider the disjoint union
\begin{equation}\label{defX}
    \wt{X}=\coprod_{n\in\N_{0}} I_{n} = I\times\N_{0}
\end{equation}
and the equivalence relation $\sim$ on $\wt{X}$ generated by the following for $n\geq1$:
\begin{align*}
	(0,n)&\sim(0,0)	\\
	(1,n)&\sim(1/n,0)
\end{align*}
Define the quotient space $\tn{HR}:=\wt{X}/\hspace{-4pt}\sim$ (see Figure~\ref{ropes}); we dub this space the \emph{Hawaiian ropes}.
\begin{figure}[h!]
    \centerline{\includegraphics{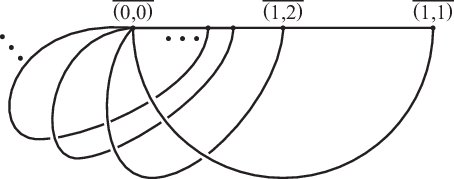}}
    \caption{Hawaiian ropes $\tn{HR}$.}
    \label{ropes}
\end{figure}
\tn{HR} is noncompact and is not a subspace of $\mathbb{R}^{3}$. Intuitively, the attached arcs are large and their interiors do not accumulate.\ss

Let $A\subset\tn{HR}$ be the image of $I\times\N$ under the quotient map $\wt{X}\to\tn{HR}$.
In other words, $A$ is the closed subspace of $\tn{HR}$ consisting of the curved arcs (`ropes') depicted in Figure~\ref{AsubX}.
\begin{figure}[h!]
    \centerline{\includegraphics{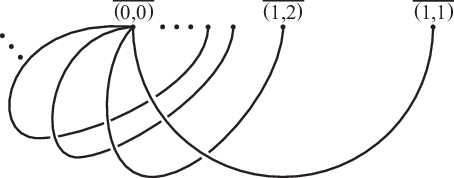}}
    \caption{Closed subspace $A\subset\tn{HR}$.}
    \label{AsubX}
\end{figure}

Our example is the quotient map
\begin{equation}\label{def_slsc_ex}
	f: \tn{HR} \to \tn{HR}/A
\end{equation}
as we now outline (see~\cite{cgm} for details).\ss

The space $\tn{HR}$ is homotopy equivalent to the wedge of countably many circles, so its fundamental group is free and countable. The quotient space $\tn{HR}/A$ is homeomorphic to the well-known Hawaiian earring (see Figure~\ref{earring}).
\begin{figure}[h!]
    \centerline{\includegraphics{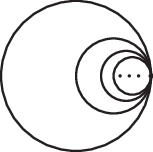}}
    \caption{Hawaiian earring $\tn{HE}\subset\mathbb{R}^{2}$.}
    \label{earring}
\end{figure}
The Hawaiian earring is arguably the simplest space that is not semilocally simply-connected; its fundamental group is uncountable (and not free)~\cite{desmit}. Thus, $f_{\sharp}$ is not surjective.\ss

Note that $\tn{HR}$ is Hausdorff, locally contractible, paracompact, and normal, but is not first countable nor metrizable~\cite[p.~12]{cgm}. The subspace $A\subset\tn{HR}$ is weakly contractible, but is not contractible nor locally path-connected.
\end{example}

\begin{remark}
In Example~\ref{Yslsc}, it is straightforward to check that the $f$-saturation of each closed subset of $\tn{HR}$ is closed in $\tn{HR}$, so $f$ is a closed map. Thus, Example~\ref{Yslsc} demonstrates the necessity of the hypothesis `fibers are locally path-connected' in the generalizations of Smale's mapping theorem in~\cite[Thm.~2]{kozlowski} and~\cite[Thm.~6.1]{dugundji}. Indeed, Example~\ref{Yslsc} satisfies all of the hypotheses of these theorems except that a single fiber, namely $A$, fails to be locally path-connected at just one point.
\end{remark}

\begin{remark}\label{Yslsc_extend}
Example~\ref{Yslsc} raises the question of whether nontriviality of $\pi_{1}(X,x_0)$ is necessary in order to create new $\pi_{1}$ in the quotient. This condition is not necessary as we now show (see~\cite{cgm} for details). Let $Z:=\tn{HR}\times I / \tn{HR}\times \{0\}$ be the cone on $\tn{HR}$.
Write $A'$ for the image of $A\times\{1\}$ in $Z$, which is a copy of $A$.
Our example is the quotient map
\[
	g: Z \to Z/ A'.
\]
The space $Z$ is contractible, locally contractible, Hausdorff, and normal, but is not metrizable. The subspace $A'\subset Z$ is closed in $Z$ and is weakly contractible. By a Mayer-Vietoris argument, the singular homology group $H_{1}(Z/A';\Z)$ is uncountable. Thus, the fundamental group of $Z/A'$ is uncountable and $g_{\sharp}$ is not surjective.
\end{remark}

\begin{example}\label{instances}
We present instances of Theorem~\ref{mainthm} where not all fibers are path-connected.
Define $h:\R\to\R$ by
\[
	h(x):= \left\{ \begin{array}
					{c@{\quad \tn{for} \quad}l}
					0	&	x\le0\\
					\sin{\frac{1}{x}}	&	x>0
					\end{array} \right.
\]
For each $c\in\R$, let $\Gamma_c$ denote the graph of $h+c$. Then $\{\Gamma_c\mid c\in\R\}$ is a partition of $\R^2$. Let $Y$ denote the associated decomposition space and let $f:(\R^2,0)\to (Y,y_0)$ be the associated quotient map. As $Y$ is indiscrete, $f$ is an instance of Theorem~\ref{mainthm} with all fibers connected but not path-connected. For examples with nontrivial $\pi_1$, consider the wedge sum $f\vee g$ where $g:(W,w_0)\to(Z,z_0)$ is, say, any map of CW-complexes that satisfies the hypotheses of Theorem~\ref{mainthm} and where $\pi_{1}(Z,z_0)\neq1$.
\end{example}

\begin{example}\label{hopf}
This example shows that nontrivial closed manifolds arise as orbit spaces of smooth (polynomial even) vector fields on euclidean space. It will be reused in Example~\ref{higherhomotopy} below.\ss

Recall that if $M$ is a smooth manifold (manifolds are Hausdorff and second-countable unless stated otherwise) and $v:M\to TM$ is a smooth vector field on $M$, then integrating $v$ yields a partition of $M$ into path-connected orbits. Each orbit is an injective image of a point, an open interval, or the unit circle. Let $\sim$ be the equivalence relation on $M$ whose elements are the orbits of $v$. The quotient space $\M$ is the \emph{orbit space} of $v$ and $q:M\to \M$ is the associated quotient map.\ss

Multiplication by $i\in\mathbb{C}$ on $\mathbb{C}^{n}$ corresponds (in real coordinates) to the following vector field on $\R^{2n}$
\[
    v(x)=\left(-x_{2},x_{1},-x_{4},x_{3},\ldots,-x_{2n},x_{2n-1}\right).
\]
Let $S^{2n-1}\subset\mathbb{R}^{2n}$ denote the unit sphere. The orbits of $v$ on $S^{2n-1}$ are great circles, one for each complex line through $0\in\mathbb{C}^{n}$. As is well-known, the orbit space of $v|_{S^{2n-1}}$ is $\CP \approx \mathbb{C}P^{n-1}$.\ss

Stereographic projection is the diffeomorphism
\[
	s: S^{2n-1} - \{p\} \to \R^{2n-1}
\]
given by
\[
    s(x)=\frac{1}{1-x_{2n}}(x_{1},x_{2},\ldots,x_{2n-1},0)
\]
where $p=(0,0,\ldots,0,1)$ and $\mathbb{R}^{2n-1}=\left\{y\in\mathbb{R}^{2n}\mid y_{2n}=0\right\}$.\ss

As $s$ is a diffeomorphism, we may push $v$ forward to obtain the vector field $u$ on $\mathbb{R}^{2n-1}$ given by
\[
    u(y)=ds_{s^{-1}(y)}(v(s^{-1}(y))).
\]
A straightforward calculation (see~\cite[pp.~25--26]{cgm}) shows that $u$ is a quadratic polynomial vector field on $\R^{2n-1}$.
One orbit of $u$ is noncompact (a properly embedded open interval) and the rest are smooth circles.\ss

Let $\R^{2n-1}/\hspace{-4pt}\sim$ denote the orbit space of the vector field $u$. We will show that $\R^{2n-1}/\hspace{-4pt}\sim$ is homeomorphic to $\mathbb{C}P^{n-1}$.\ss

Consider the diagram
\[
\xymatrix{
        \mathbb{R}^{2n-1}   \ar[r]^-{s^{-1}}    \ar[d]_{q'} &   S^{2n-1}    \ar[d]^{q}\\
        \mathbb{R}^{2n-1}/\hspace{-4pt}\sim  \ar@{-->}[r]^-{f}   & \CP}
\]
where $q$ and $q'$ are the associated quotient maps. By construction, $q\circ s^{-1}$ is constant on each fiber of $q'$, so the universal property of quotient maps implies that there exists a unique continuous function $f$ making the diagram commute. Plainly, $f$ is a bijection. Let $D=\left\{x\in S^{2n-1}\mid x_{2n}\le 0\right\}$ be the lower hemisphere of $S^{2n-1}$. Notice that each orbit in $S^{2n-1}$ intersects $D$ in at least a nontrivial arc of points and that $D$ is compact. Thus $s(D)$ is compact and maps by $q'$ surjectively to $\mathbb{R}^{2n-1}/\hspace{-4pt}\sim$. Thus $\mathbb{R}^{2n-1}/\hspace{-4pt}\sim$ is compact. As $\CP$ is Hausdorff, $f$ is a homeomorphism as desired. Using induced functional structures (see~\cite[pp.~69--72]{bredon_topology}), the fiber bundle structures permit one to show that $f$ is actually a diffeomorphism of smooth manifolds.
\end{example}

\begin{example}\label{higherhomotopy}
We construct an instance of Theorem~\ref{mainthm} whose domain is contractible and locally contractible, whose fibers are contractible, and whose quotient, while necessarily simply-connected, has higher homotopy and homology.
This shows that the analogy of our Theorem~\ref{mainthm} with the long exact sequence of a fibration does not extend to higher homotopy groups, so Theorem~\ref{mainthm} at this level of generality is specific to the fundamental group.\ss

Fix notation as in Example~\ref{hopf} with $n=2$, so $q$ is the classical Hopf fibration. Let $\sim$ denote the associated equivalence relation on $S^{3}$. The unique complex line contained in the set $\left\{(x_{1}+ix_{2},x_{3}+ix_{4})\mid x_{4}=0\right\}$ is $\mathbb{C}\times\left\{0\right\}$. Thus $D$ intersects orbits of $q$ in the following way: one intersection is a complete orbit $\left\{(x_{1},x_{2},0,0)\mid x_{1}^{2}+x_{2}^{2}=1\right\}$ and the rest are semicircles each intersecting $\partial{D}$ in a pair of antipodal points. Let $D'=D-\left\{(1,0,0,0)\right\}$ and let $\sim$ denote the restriction of $\sim$ to $D'$ where no confusion should arise. We have the diagram
\[
\xymatrix{
        D'  \ar[r]^-{j} \ar[d]_{\pi}    &   S^{3}   \ar[d]^{q}\\
        D'/\hspace{-4pt}\sim \ar@{-->}[r]^-{g}   &   S^{3}/S^{1}}
\]
where $j$ is inclusion and $\pi$ is a quotient map. The composition $q\circ j$ is constant on each fiber of $\pi$ and so the universal property of quotient maps implies that the unique function $g$ making the diagram commute is continuous. It is easy to see that $g$ is a bijection. We claim that $D'/\hspace{-4pt}\sim$ is compact. To see this, let $B\subset\mathbb{R}^{4}$ denote the open ball of radius $1/2$ centered at $(1,0,0,0)$. Note that $D'-B$ is a compact subset of $D'$. Note also that each fiber of $\pi$ has diameter $2$ and thus intersects $D'-B$ nontrivially. Therefore $\left.\pi\right|_{D'-B}$ is surjective and $D'/\hspace{-4pt}\sim$ is compact. As $S^{3}/S^{1}$ is Hausdorff, $g$ is a homeomorphism. Thus we have $D'$, which is diffeomorphic to closed upper $3$--space, partitioned into connected arcs (one open and the rest closed) with quotient $S^{2}$, completing the example. Note that stereographic projection from $(1,0,0,0)$ explicitly exhibits a diffeomorphism from $D'$ to half-space, with the partition induced by a quadratic vector field (cf. Example~\ref{hopf}), so this example is remarkably simple from an analytic viewpoint.
\end{example}

\section{Proof of Theorem~\ref{mainthm}}\label{sec:proof}

By hypothesis, $f:\pa{X,x_{0}}\to\pa{Y,y_{0}}$ is a quotient map with connected fibers, $X$ is locally path-connected, and $Y$ is semilocally simply-connected. We must show that $f_{\sharp}:\pi_{1}\pa{X,x_{0}}\to\pi_{1}\pa{Y,y_{0}}$ is surjective.\ss

First, $Y$ is locally path-connected since it is a quotient space of a locally path-connected space. Second, restricting $f$ to the path components of $X$ and $Y$ containing the basepoints $x_0$ and $y_0$, we can and do further assume $X$ and $Y$ are path-connected.
Indeed, each component $C$ of $X$ is open and closed in $X$, and is $f$-saturated (since fibers are connected). Hence, $f(C)$ is open and closed in $Y$, $f(C)$ is a path component of $Y$, and $f|_{C}:C\to f(C)$ is a quotient map.\ss

Now, we apply classical covering space theory to $Y$, which requires no Hausdorff hypothesis~\cite[\S~1.3]{hatcher}. Corresponding to the subgroup $\tn{Im}(f_{\sharp})$ of $\pi_1 \pa{Y,y_0}$, there exists a based covering map
\hbox{$c:\pa{Z,z_0}\to\pa{Y,y_0}$},
unique up to isomorphism, such that $Z$ is connected and
\begin{equation}\label{groups}
	\tn{Im} (c_{\sharp}) = \tn{Im} (f_{\sharp}).
\end{equation}
Equation~\eqref{groups} yields a (unique) lift $g$ of $f$, so the following diagram commutes.
\begin{equation}\begin{split}\label{lift}
\xymatrix{
    &   \pa{Z,z_0}   \ar[d]^-{c}\\
    \pa{X,x_0}   \ar[r]^-{f}   \ar@{-->}[ur]^-{g}  &   \pa{Y,y_0}}
\end{split}\end{equation}
We will show that the restriction
\begin{equation}\label{rest}
	c|_{\tn{Im}(g)} : \tn{Im}(g) \to Y
\end{equation}
is a homeomorphism onto $Y$.\ss

Surjectivity of~\eqref{rest} follows by commutativity of~\eqref{lift} since $f$ is surjective. For injectivity, note that fibers of $c$ are discrete and fibers of $f$ are connected. Commutativity of~\eqref{lift} implies that $g$ maps each fiber of $f$ into a fiber of $c$. Hence, $g$ is constant on each fiber of $f$ and~\eqref{rest} is injective.\ss

We claim that $\tn{Im}(g)$ is open in $Z$. Let $g(x)\in \tn{Im}(g)$. Then $f(x)\in Y$ lies in a connected open set $U$ evenly covered by $c$. Note that $c^{-1}(U)$ is the disjoint union
\[
    c^{-1}(U)=\coprod_{j\in J}U_{j}
\]
of connected open sets in $Z$ where $J$ is some index set, $c|_{U_{j}}:U_{j}\rightarrow U$ is a homeomorphism for each $j\in J$, and $g(x)\in U_{0}$. The connected components of $f^{-1}(U)$ are $f$-saturated and open in $X$, so their images under $f$ are disjoint and open. Therefore $f^{-1}(U)$ is connected, $f$-saturated, and open in $X$. In particular, $g$ maps $f^{-1}(U)$ into $U_{0}$. To see that $g$ maps $f^{-1}(U)$ onto $U_{0}$, let $z\in U_{0}$. Then $c(z)\in U$ and so there exists $p\in f^{-1}(U)$ such that $f(p)=c(z)$. Now $g(p)\in U_{0}$ and commutativity of~\eqref{lift} implies $c(g(p))=c(z)$. But $c|_{U_{0}}:U_{0}\to U$ is a homeomorphism and so $g(p)=z$ as desired. Hence $g(f^{-1}(U))=U_{0}$, which is an open neighborhood of $g(x)$ in $Z$, and $\tn{Im}(g)$ is open in $Z$, proving the claim.\ss

As $c$ is a local homeomorphism and $\tn{Im}(g)$ is open in $Z$, \eqref{rest} is a local homeomorphism. Thus~\eqref{rest}, being a bijective local homeomorphism, is a homeomorphism as desired.\ss

At the level of based loops, the homeomorphism~\eqref{rest} implies that $c_{\sharp}$ is surjective. Equation~\eqref{groups} implies that $f_{\sharp}$ is surjective, completing the proof of Theorem~\ref{mainthm}. \qed \medskip

\begin{remark}
As $c$ is a based covering map, $c_{\sharp}$ is injective. Coupled with the proof above, we see that $c_{\sharp}$ is an isomorphism and so $c$ is a one sheeted cover. By a based covering space isomorphism, we may take $\pa{Z,z_0}=\pa{Y,y_0}$ and $c=\tn{id}$. By commutativity of~\eqref{lift}, $g=f$ and so $g$ is surjective. The proof of Theorem~\ref{mainthm} may be reorganized to argue directly that $g$ is surjective.
\end{remark}

\section{Theorem~\ref{mainthm} and the Simplicial Category}\label{simplicial_category}

On one hand, Theorem~\ref{mainthm} already applies to any surjective simplicial map with connected fibers (see Corollary~\ref{simplicial_surjection_quotient} below). On the other hand, for simplicial maps one may proceed using purely simplicial techniques and achieve somewhat more (see Corollary~\ref{simplicial_mainthm} below). The reader may contrast Corollary~\ref{simplicial_mainthm} with results of Quillen~\cite[Props.~1.6,~7.6]{quillen}. Under much stronger hypotheses, Quillen proves stronger conclusions. Our Corollary~\ref{simplicial_mainthm} is more in the spirit of Theorem~\ref{mainthm} and results of Smale and others from Section~\ref{s:intro}. \ss

Simplicial complexes (see~\cite[Ch.~3]{spanier}) will be unordered, and not necessarily finite or even locally finite.
Simplices are closed unless stated otherwise.
The topological space underlying the simplicial complex $K$, and endowed with the coherent topology, is denoted $\left|K\right|$.
The \hbox{$n$-skeleton} of $K$ is denoted $K^{(n)}$.
If $x\in\left|K\right|$, then the \emph{carrier} of $x$ in $K$, denoted $\tn{carr}(x,K)$, is the unique smallest simplex $\sigma\in K$ such that $x\in\left|\sigma\right|$.
If $x$ is the barycenter of a simplex $\sigma\in K$, then $\tn{carr}(x,K)=\sigma$.
Given a simplicial map $f:K\to L$, there is an associated map of topological spaces $\left|f\right|:\left|K\right|\to\left|L\right|$.

\begin{lemma}\label{n-simplex}
Let $f:K\to L$ be a simplicial map and let $\tau\in L\cap\tn{Im}(f)$ be an $n$-simplex. Then there exists an $n$-simplex $\sigma\in K$ such that $f|_{\sigma}:\sigma\to\tau$ is a bijection.
\end{lemma}

\begin{proof}
Let $y$ be the barycenter of $\tau$. Let $x\in \ts{f}^{-1}(y)$ and let $\mu=\tn{carr}(x,K)$. As $f(\mu)\in L$ and $y\in\ts{f(\mu)}$, we see that $\tau=\tn{carr}(y,L)$ is a face of $f(\mu)$. Some face $\sigma$ of $\mu$ maps bijectively by $f$ to $\tau$.
\end{proof}

\begin{corollary}\label{simplicial_surjection_quotient}
If $f:K\to L$ is a surjective simplicial map, then $\ts{f}$ is a quotient map.
\end{corollary}

\begin{proof}
Let $V\subset\left|L\right|$ such that $\ts{f}^{-1}(V)$ is open in $\left|K\right|$. Let $\tau\in L$. We must show that $\ts{\tau}\cap V$ is open in $\ts{\tau}$. So, let $y\in\ts{\tau}\cap V$. By Lemma~\ref{n-simplex}, there is $\sigma\in K$ such that $f|_{\sigma}:\sigma\to\tau$ is a bijection. As $\ts{\sigma}\cap \ts{f}^{-1}(V)$ is open in $\ts{\sigma}$ and $\ts{f}|_{\ts{\sigma}}$ is a homeomorphism, $\ts{f}(\ts{\sigma}\cap \ts{f}^{-1}(V))\subset V$ is an open neighborhood of $y$ in $\ts{\tau}$. 
\end{proof}

\begin{remark}\label{not_open}
While the geometric realization $\ts{f}$ of a simplicial surjection $f$ is a quotient map, it is not necessarily open or closed.
Let $L=\{ \{v_1\}, \{v_2\}, \{v_1,v_2\} \}$ and let
\[
	K=\{ \{u_0\}, \{u_1\}, \{u_2\}, \{u_0,u_1\}, \{u_1,u_2\} \}.
\]
Define $f:K\to L$ by $f(u_0)=v_1$, $f(u_1)=v_1$, and $f(u_2)=v_2$.
Let $\sigma=\{u_0,u_1\}$.
The open simplex $\left\langle \sigma \right\rangle$ does not have open image in $\ts{L}$.
For a non-closed example, let
\[
	M=\{ \{x_1\}, \{x_2\}, \{x_3\},\ldots, \{x_1,x_2\}, \{x_2,x_3\}, \ldots \}.
\]
Define $g:M\to L$ by $g(x_{2i-1})=v_1$ and $g(x_{2i})=v_2$ for $i\in\N$.
The discrete set
\[
	\left\{ \left. \frac{1}{i+1}x_{i}+\frac{i}{i+1}x_{i+1} \right| i\in\N \right\}
\]
is closed in $\ts{M}$ but its image is not closed in $\ts{L}$.
However, if $f$ is a simplicial map of locally finite complexes and $\ts{f}$ is a proper map, then $\ts{f}$ is a closed map (cf. the paragraph following Theorem~\ref{mainthm} in Section~\ref{s:intro}).
\end{remark}

Next, we recall a few definitions (cf.~\cite{smale_a}). Let $I=[0,1]$.
A map $f:X\to Y$ has the \emph{covering homotopy property for a point}, abbreviated CHP,
provided: for each map $\alpha:I\to Y$ and $x\in f^{-1}(\alpha(0))$, there exists a map
$\beta:I\to X$ such that $\beta(0)=x$ and $f\circ\beta=\alpha$.
One says $f$ satisfies the \emph{covering homotopy property for a point up to homotopy},
abbreviated CHPH, provided: given $f$, $\alpha$, and $x$ as above, there exists $\beta:I\to X$
such that $\beta(0)=x$, $f(\beta(1))=\alpha(1)$, and $f\circ\beta$ is path homotopic to $\alpha$.\ss

The CHP fails to hold even for the very simple map $\ts{f}$ in Remark~\ref{not_open} above (this map is closed, surjective, and has contractible fibers): an injective path from $v_1$ to $v_2$ has no lift beginning at $u_0$. But, the CHPH is useful for simplicial maps. The following lemma is a simplicial CHPH.

\begin{lemma}\label{simplicial_chph}
Let $f:K\to L$ be a simplicial map such that $L^{(1)}\subset\tn{Im}(f)$ and $f^{-1}(v)$ is connected for each vertex $v\in L^{(0)}$. Let $\alpha:I\to\ts{L}$ be induced by a simplicial map. Choose arbitrary vertices $v_{0},v_{1}\in K^{(0)}$ such that $f(v_0)=\alpha(0)$ and $f(v_1)=\alpha(1)$. Then there exists a map $\beta:I\to\ts{K}$, induced by a simplicial map, such that $\beta(0)=v_0$, $\beta(1)=v_1$, and $\ts{f}\circ\beta$ is path homotopic to $\alpha$.
\end{lemma}

\begin{proof}
To say $\alpha$ is induced by a simplicial map means there exists a triangulation $t:\ts{J}\to I$ and a simplicial map $a:J\to L$ such that $\alpha=\ts{a}\circ t^{-1}$. Without loss of generality, we assume $a$ is not constant on any edge of $J$.\ss

Let $J^{(0)}=\{w_0,w_1,\ldots,w_n\}$ where $a(w_0)=\alpha(0)$ and $a(w_n)=\alpha(1)$, and let $\{e_1,e_2,\ldots,e_n\}$ be the edges of $J$ where $e_i=\{w_{i-1},w_i\}$. Then
\[
	\tau_i :=a(e_i)=\{a(w_{i-1}),a(w_i)\}
\]
is an edge of $L$. As $L^{(1)}\subset\tn{Im}(f)$, Lemma~\ref{n-simplex} yields an edge
$\sigma_i =\{\sigma_i^{-},\sigma_i^{+}\}$
of $K$ such that $f(\sigma_i^{-})=a(w_{i-1})$ and $f(\sigma_i^{+})=a(w_i)$. Consider adjacent edges $e_i$ and $e_{i+1}$ of $J$. Notice that $\sigma_i^{+}$ and $\sigma_{i+1}^{-}$ lie in $f^{-1}(a(w_i))$, a connected subcomplex of $K$. Thus, we may link $\sigma_i^{+}$ and $\sigma_{i+1}^{-}$ by a simplicial path $C_i$ in $f^{-1}(a(w_i))$. Similarly, $v_0$ and $\sigma_1^{-}$ lie in $f^{-1}(a(w_0))$ and $v_1$ and $\sigma_n^{+}$ lie in $f^{-1}(a(w_n))$, yielding the linking paths $C_0$ and $C_n$ respectively. Then $C_0 \sigma_1 C_1 \sigma_2 C_2 \cdots \sigma_n C_n$ is a simplicial path as desired.
\end{proof}

Lemma~\ref{simplicial_chph} and the Simplicial Approximation theorem immediately imply the following simplicial analogue of Theorem~\ref{mainthm}.

\begin{corollary}\label{simplicial_mainthm}
If $f:(K,k_0)\to(L,l_0)$ is a simplicial map such that $L^{(1)}\subset\tn{Im}(f)$ and $f^{-1}(v)$ is connected for each vertex $v\in L^{(0)}$, then the induced homomorphism $f_{\sharp}$ is surjective.
\end{corollary}

\section{Concluding Remarks}\label{conclusion}

The opening question in Arnol{\textprime}d's \emph{Problems on singularities and dynamical systems}~\cite[p.~251]{arnold} asks whether an exotic $\R^4$ may appear as the orbit space of a polynomial or trigonometric vector field on $\R^5$.
If one asks merely for a smooth vector field, then the answer is affirmative for every exotic $\R^4$ as noted by Arnol{\textprime}d. Thus the gist of the question is: \emph{does dynamics produce exotic differentiable manifolds in the simplest possible scenario?} This harks back to the quadratic polynomial vector field on $\R^3$ producing the Lorenz attractor: chaos is exhibited by a continuous dynamical system in the simplest possible setting.\ss

Inspired by Arnol{\textprime}d's question, it is natural to ask: for a given manifold $M$, which manifolds or spaces may appear as orbit spaces $\M$ of vector fields on $M$? Of course, orbit spaces are often non-Hausdorff. Still, nontrivial closed manifolds arise even for quadratic polynomial vector fields on euclidean space (see Example~\ref{hopf} above). The case where $M$ is the euclidean plane has been studied in, e.g.,~\cite{haefliger_reeb} and~\cite{gabai_kazez}.\ss

For a given manifold $M$, Theorem~\ref{mainthm} gives an obstruction to spaces arising as orbit spaces $\M$. Manifolds arising as orbit spaces of $\R^n$, for example, must be simply-connected. Theorem~\ref{mainthm} also implies that the line with two origins~\cite[p.~15]{cg}, a well-known non-Hausdorff smooth $1$-manifold with fundamental group isomorphic to $\Z$, is not an orbit space of $\R^2$. Theorem~\ref{mainthm} may also be applied to $p$-dimensional foliations of $M$ and their associated leaf spaces.\ss

The question arises whether every orbit space $\M$ of a smooth vector field is semilocally simply-connected. The authors know of no counterexample. In~\cite{cg} the first two authors give an affirmative answer for a large class of gradient vector fields of Morse functions, along with some related results.\ss

In general, orbit space quotient maps $q:M\to \M$ do not have the CHP (see Section~\ref{simplicial_category} above for definitions and see Example~4.6 of~\cite{cg}).
We ask whether all orbit space quotient maps have the CHPH.\ss

From a dynamical systems viewpoint, orbit spaces are important, of course, since they are topological invariants of continuous flows. It seems intriguing to ask about the particular dynamical significance of our local/global fundamental group questions.

\end{document}